\documentclass[12pt]{article}

\usepackage{amssymb}

\newcommand{\beqn}{\begin{eqnarray}}
\newcommand{\eeqn}{\end{eqnarray}}
\newcommand{\be}{\begin{equation}}
\newcommand{\ee}{\end{equation}}
\newcommand{\ba}{\begin{array}}
\newcommand{\ea}{\end{array}}

\newcommand{\re}{\ref}
\newcommand{\ci}{\cite}
\newcommand{\la}{\label}

\newcommand{\ds}{\displaystyle}
\newcommand{\ve}{\varepsilon}

\begin{document}

\renewcommand{\theequation}{\thesection.\arabic{equation}}
\newtheorem{theorem}{Theorem}[section]
\renewcommand{\thetheorem}{\arabic{section}.\arabic{theorem}}
\newtheorem{definition}[theorem]{Definition}
\newtheorem{deflem}[theorem]{Definition and Lemma}
\newtheorem{lemma}[theorem]{Lemma}
\newtheorem{example}[theorem]{Example}
\newtheorem{remark}[theorem]{Remark}
\newtheorem{remarks}[theorem]{Remarks}
\newtheorem{cor}[theorem]{Corollary}
\newtheorem{pro}[theorem]{Proposition}
\mathsurround=2pt

\newcommand{\R}{{\mathbb{R}}}
\newcommand{\C}{{\mathbb{C}}}
\newcommand{\Z}{{\mathbb{Z}}}
\newcommand{\N}{{\mathbb{N}}}
\newcommand{\pr}{\prime}

\begin{center}
{\Large\bf Boundedness and Convergence of Solutions \\
\medskip
for the String Coupled to a Nonlinear Oscillator}
\vspace{1cm}

{\large T.V.~Dudnikova}\\
{\it Keldysh Institute of Applied Mathematics RAS\\
 Moscow 125047, Russia}\\
e-mail:~tdudnikov@mail.ru
\end{center}
\vspace{1cm}

 \begin{abstract}
A system of equations consisting of an infinite string
coupled to a nonlinear oscillator is considered.
The Cauchy problem for the system
with the periodic initial data is studied.
The main goal is to prove the convergence of the solutions
as $t\to\infty$ to a time periodic solution.
\bigskip\\
{\it Key words and phrases}:
an infinite string coupled to a nonlinear oscillator,
the Cauchy problem, periodic initial data,
the limit amplitude principle
\end{abstract}

\section{Introduction}

Consider the following problem for a function $u(x,t)\in C(\R^{2})$:
\beqn\label{1}
(\mu+m\delta(x))\ddot u(x,t)=
\kappa u''(x,t)+\delta(x)F(u(x,t)),\quad t\in\R,\quad x\in \R.
\eeqn
Here $m\ge0$,
$\mu,\kappa>0$;  $\dot u\equiv{\partial u/\partial t}$,
$ u'\equiv \partial u/\partial x$.
The initial data (when $t=0$) for Eq.~(\ref{1}) are assumed to be
periodic, see Definition \ref{d1.5} below.

By definition, Eq.~(\ref{1}) is equivalent to the following system:
\beqn
\mu\ddot u(x,t)=\kappa u''(x,t),\quad t\in\R,
\quad x\in \R\setminus \{0\},\la{3}\\
m\ddot y(t)=F(y(t))+\kappa\,[u'(0+,t)-u'(0-,t)],\quad t\in\R, \label{4}
\eeqn
where
\be\label{4'}
y(t)=u(0-,t)=u(0+,t),\quad t\in\R.
\ee
Physically, the system describes small crosswise
oscillations of an
infinite string stretched parallel to the $Ox$-axis.
$\mu$ is the line density of the string, $\kappa$ is its tension,
$F(y)$ is an external (nonlinear, in general) force field perpendicular to $Ox$.
In the case $m=0$, the string is coupled
 to a  spring of a rigidity $F(y)$.
In the case $m>0$, a ball of mass $m$ is attached to the string
at the point $x=0$, and the field $F(y)$ subjects the ball.

The system (\ref{3})--(\ref{4'}) was considered first  by Lamb~\cite{La}
for the linear case, i.e., when $F(y)=-ry$ with a positive constant $r$.
For general nonlinear functions $F(y)$, this model
was studied  by Komech in the paper \cite{KMo},
where the transitions to stationary states were established for
finite energy solutions.
In the present paper, the solutions of infinite energy with
space-periodic initial data are considered.
Main goal is to prove  that each solution $u(x,t)$ to the system
for large times is close to a time-periodic solution
(see Theorem~\ref{t1} below).
\medskip

 Let us describe our assumptions on the external force $F(y)$.

Denote by $V(y)= -\ds\int F(y)\,dy$ the potential energy of the external
field, $F(y)= -V'(y)$, $y\in \R$.
We assume that
\beqn\label{F}
 F(y)\in C^1(\R),\quad F(y)\to\mp\infty
\quad\mbox {as } y\to\pm\infty.
\eeqn
Obviously, condition (\ref{F}) implies that
\beqn\label{V}
V\in C^2(\R),\quad
 V(y)\to \infty \quad \mbox {as }\, |y|\to\infty.
\eeqn

Let us introduce a class ${\cal E}$ of solutions $u(x,t)$
to Eq.~(\ref{1}) with locally finite energy.
\begin{definition}\la{d2}
A function $u(x,t)$ belongs to ${\cal E}$
if $u\in C(\R^2)$ and $\dot u, u'\in L^2_{{\rm loc}}(\R^2)$,
where the derivatives are understood in the sense
of distributions.
\end{definition}

For $u(x,t)\in{\cal E}$, the system (\ref{3})--(\ref{4})
is understood as follows (see \cite{KMo}).

For $u\in C(\R^2)$, Eq.~(\re{3}) is understood in the sense
of distributions in the region $(x,t)\in \R^2$, $x\ne 0$.
Moreover,  Eq.~(\re{3}) is
equivalent to the d'Alembert decomposition ($a=\sqrt{\kappa/\mu}$)
\be\la{5}
u(x,t)=f_\pm(x-at)+ g_\pm(x+at),\quad \pm x>0,\quad t\in\R,
\ee
where $f_{\pm},g_{\pm}\in C(\R)$, since  $u(x,t)\in C(\R^2)$.

We now explain Eq.~(\re{4}).  Equality (\re{5}) implies
$$
u'(x,t)=f'_\pm(x-at)+g'_\pm(x+at),\quad
 \pm x>0,\quad t\in\R,
$$
where all derivatives are understood in the sense of distributions.
For $u(x,t)\in C(\R^2)$ satisfying (\re{3}), write
\be\la{6}
u'(0\pm,t):=f'_\pm(-at)+g'_\pm(at). 
\ee
Note that condition $u(x,t)\in{\cal E}$ implies that
$f'_\pm, g'_\pm\in L^2_{{\rm loc}}(\R)$.
The derivative $\ddot y(t)$ of $y(t)\in C(\R)$
is understood in the sense of distributions.
Moreover, for $m\not=0$, Eq.~(\ref{4})
and condition (\ref{F}) imply that
$\ddot y(t)=\ddot u(0\pm,t)\in L^2_{{\rm loc}}(\R)$. Hence,
if $m\not=0$, $y(t)\in C^1(\R)$ for any solution $u\in{\cal E}$.
\medskip

We study the Cauchy problem for the system (\ref{3})--(\ref{4})
with the initial conditions
\beqn
\label{2'}
u|_{t=0}=u_0(x),\quad \dot u|_{t=0}=u_1(x),\quad x\in\R,\\
\label{2''}
\dot y|_{t=0}=y_1\quad (\mbox{if }\, m\not=0).
\eeqn
We assume that $y_1\in\R$ and the initial data $u_0(x), u_1(x)$
belong to the space ${\cal H}$.
\begin{definition}
The pair of functions $(u_0,u_1)$ belongs to the space ${\cal H}$ if
$u_0\in C(\R)$,  $u'_0, u_1\in L^2_{{\rm loc}}(\R)$.
\end{definition}
\begin{pro}\la{p1}
Let condition (\ref{V}) hold and $(u_0,u_1)\in{\cal H}$, $y_1\in\R$.
Then the Cauchy problem (\ref{3})--(\ref{4'}), (\ref{2'}), (\ref{2''})
has a unique solution $u(x,t)\in {\cal E}$.
\end{pro}

This proposition is proved in Section \ref{sec2}.
\medskip

To prove the main result
we impose additional conditions on the initial data $(u_0,u_1)$.
At first, for an $\omega>0$, we introduce a class $P^\omega$
of the  space periodic functions.

\begin{definition}
For $\omega>0$, we say that $u\in P^\omega$ if
$u(x\pm\omega)=u(x)$ for $\pm x>0$.
\end{definition}
\begin{definition}\label{d1.5}
For $\omega>0$,
$(u_0,u_1)\in{\cal H}^\omega$ if $u_0\in C^1(\R)$, $u_1\in C(\R)$
and $u_0, u'_0, u_1\in P^\omega$.
\end{definition}

In the case $m=0$, the main result is the following convergence theorem.
\begin{theorem}\la{t1}
Let  $m=0$, condition (\ref{F}) hold
and $(u_0,u_1)\in{\cal H}^\omega$ for some $\omega>0$.
Then for every solution  $u(x,t)\in {\cal E}$
of the Cauchy problem (\ref{3})--(\ref{4'}), (\ref{2'}) there exists a
 solution $u_p(x,t)\in {\cal E}$ to Eq.~(\ref{1}) such that
\be\label{10}
u_p(x,t+\omega/a)=u_p(x,t)\quad\mbox{for }\quad (x,t)\in\R^2:\,|t|>|x|/a,
\ee
 and for every $R>0$,
\beqn\la{9}
\int\limits_{|x|<R}\!\!\left(
|\dot u(x,t)-\dot u_p(x,t)|^2+|u'(x,t)-u'_p(x,t)|^2\right)dx
+\max_{|x|<R} |u(x,t)-u_p(x,t)|\to 0
\eeqn
as $t\to\infty$.
\end{theorem}

This theorem is proved in Section~\ref{sec3}.
The similar result holds for $m\not=0$ under additional restrictions
on the function $F(y)$ (see Section \ref{sec3}).


In Appendix~B we consider Eq.~(\ref{1})  for $t>0$
under the initial condition
\be\label{2}
u(x,t)|_{t\le0}=p(x+at),\quad x\in \R,
\ee
where the function $p(z)\in P^\omega$, $p\in C^1(\R)$,
$p(x)=p_0$ for  $x\leq 0$, and $F(p_0)=0$.
In this case, the convergence (\ref{9}) holds, i.e.,
 the solution $u(x,t)$ of the problem (\ref{3})--(\ref{4'}), (\ref{2})
either is a time-periodic for $|x|\le at$ with
period $\omega/a$ or converges to a
function $u_p(x,t)\in{\cal E}$ satisfying (\ref{10}).
Moreover, the function $u_p(x,t)$ is a solution of Eq.~(\ref{1}) for $t>0$
under the condition $u_p(x,t)|_{t\le0}=q(x+at)$.
Here $q(x)=q_0$ for $x\le0$ and $q(x)=q_0+p(x)-p_0$ for $x>0$,
with some point $q_0\in\R$ depending on $p_0$.
\medskip


We outline the strategy of the proof of (\ref{9}).
At first, using the d'Alembert method,
we reduce the  problem (\ref{3})--(\ref{4'}),
(\ref{2'}), (\ref{2''})
to the study of the following Cauchy problem for the function $y(t)$,
\be\label{100}
m\ddot y +(2\kappa/a)\dot y-F(y(t))=2\kappa p'(at),\quad t\in\R,
\ee
with some $\omega$-periodic function $p$ (see formula (\ref{1.6}) below)
and with the initial conditions
\be\label{ll}
\ba{l}
y|_{t=0}=y_0=u_0(0),\\
\dot y|_{t=0}=y_1\,\, (\mbox{if }\,m\not=0).
\ea\ee
Further, for $m=0$, we show (see Theorem~\ref{t2.1})
 that any solution of Eq.~(\ref{100})
either $\omega/a$-periodic or tends to an $\omega/a$-periodic
solution $y_p(t)$,
i.e., $|y(t)-y_p(t)|\to0$ as $t\to\infty$.
Finally, using the explicit formula (\ref{1.11})
for $u(x,t)$ we derive the results
of Theorem~\ref{t1}.
\medskip

If $m\not=0$, the behavior of solutions to Eq.~(\ref{100}) is
more complex.
If $F(y)=-ax-by^3$, the equation of the form (\ref{100})
is called the {\em Duffing equation with damping}, see for example, \cite{Lo55,Sh}.
Eq.~(\ref{100}) is a particular case of the  generalized Li\'enard equations
with a forcing term $e(t)=2\kappa p'(at)$,
\be\label{1.16'}
\ddot y+f(y)\dot y+g(y)=e(t).
\ee
Eq.~(\ref{1.16'}) with $g(y)=y$ and $e(t)\equiv0$
was studied first by Li\'enard \cite{Li}.
A class of equations of the form (\ref{1.16'})
has been widely investigated in the literature, see, for example,
Cartwright~\cite{C}, 
Littlewood~\cite{CL}, 
Levinson~\cite{Lev},
Loud~\cite{Lo55, Lo57}, Reuter \cite{Re}.
We refer the reader to the survey works \cite{Le, P0, P, RSC, SC}
for a detailed discussion of the results and methods
concerning these equations.
Some results concerning Eq.~(\ref{100}) are given in Section \ref{sec3}.
In particular, condition~(\ref{V}) implies that for large times
the pairs $Y(t)=(y(t),\dot y(t))$ (where $y(t)$ is a solution
of (\ref{100})) belong to a fixed bounded region of
$\R^2$.
Denote by $U(t,0)$ the solving operator to the Cauchy problem (\ref{100}),
(\ref{ll}). By the Pliss results  \cite{P0, P}, there exists a set $I\subset \R^2$
which is invariant w.r.t. $U(\omega/a,0)$.
Moreover, the set $I$ is not empty and has zero Lebesgue measure.
Introduce an integral set ${\cal S}\subset \{(Y(t),t)\in\R^3\}$
consisting of the solutions of Eq.~(\ref{100})
with the initial values $(y_0,y_1)\in I$.
Let ${\cal S}_\tau$ denote the intersection
of ${\cal S}$ and the hyperplane $t=\tau$,
and $\rho(Y, {\cal S}_\tau)$ stand for the distance
between a point $Y\in\R^2$ and the set ${\cal S}_\tau$.
In Section~\ref{sec3} we check that every solution of Eq.~(\ref{100})
tends to the set ${\cal S}$ as $t\to\infty$, i.e.,
$\rho(Y(\tau),{\cal S}_\tau)\to0$ as $\tau\to\infty$.
Hence the explicit formula (\ref{1.11}) for the solutions $u(x,t)$
implies that for any $R>0$,
\beqn\la{11}
\inf\Big\{\int\limits_{|x|<R}\!\! \left(
|\dot u(x,t)-\dot u_p(x,t)|^2+|u'(x,t)-u'_p(x,t)|^2\right)\,dx
+\max_{|x|<R} |u(x,t)-u_p(x,t)|\Big\}
\eeqn
vanishes as $t\to\infty$, where the infinitum is taken over all
solutions $u_p(x,t)\in{\cal E}$ of the problem (\ref{3})--(\ref{4})
such that $u_p(0\pm,t)=y_p(t)$ and $(y_p(t), \dot y_p(t))\in {\cal S}_t$.

We give additional restrictions on the function $F(y)$ (see Examples \ref{ex1}--\ref{ex3})
when the set $I$ has a unique point
and then Eq.~(\ref{100}) has a unique stable periodic solution. In this case,
every solution of Eq.~(\ref{100}) tends to a $\omega/a$-
periodic solution $y_p(t)$ as $t\to\infty$, and convergence (\ref{9}) holds.

\setcounter{equation}{0}
\section{Existence of solutions}\label{sec2}

In this section we prove Proposition \re{p1}.
The method of  construction of finite energy solutions to the Cauchy problem
(\ref{3})--(\ref{4'}), (\ref{2'}), (\ref{2''}) was given by Komech in \ci{KMo}.
We apply this method to the infinite energy solutions.
For simplicity, we consider only the case $t>0$.
Substituting (\re{5}) into initial conditions (\re{2'}), we have
\beqn\la{1.1}
\ba{ll}
f_\pm(z)=u_0(z)/2-1/(2a)\ds\int_0^z u_1(y)\,dy+C_\pm,&\mbox{for }\pm z>0,\\
g_\pm(z)=u_0(z)/2 +1/(2a)\ds\int_0^z u_1(y)\,dy-C_\pm,& \mbox{for }\pm z>0,
\ea
\eeqn
where we can put constants $C_\pm=0$.
On the other hand, substituting (\ref{5}) into the condition~(\ref{4'}), we have
\beqn\la{1.2}
y(t)=f_-(-at)+g_-(at)=f_+(-at)+g_+(at)\quad\mbox{for }\,t\in\R.
\eeqn
By (\ref{1.2}), we can determinate $g_-(z)$ with $z>0$  and $f_+(z)$ with $z<0$ as follows:
\beqn\la{1.3}
g_-(z)=y(z/a)-f_-(-z),\quad f_+(-z)=y(z/a)-g_+(z)\quad \mbox{for }\,z>0.
\eeqn
Therefore, for $t>0$ we obtain
\beqn\label{1.11}
u(x,t)=\left\{\ba{lll}
f_+(x-at)+g_+(x+at)& \mbox{for }\,\,x\ge at\\
y(t-x/a)+g_+(x+at)-g_+(at-x)& \mbox{for }\,\,0\le x<at\\
y(t+x/a)+f_-(x-at)-f_-(-at-x)&\mbox{for }\,\, -at\le x<0\\
f_-(x-at)+g_-(x+at)& \mbox{for }\,\,x<-at
\ea
\right.
\eeqn
where $f_\pm\in C(\R_\pm)$, $f'_\pm\in L^2_{{\rm loc}}(\R_\pm)$
with $\R_\pm=\{x\in\R:\pm x>0\}$.
Moreover, by definition~(\ref{6}), we have
\beqn\nonumber
\ba{rcl}
u'(0+,t)&:=&f'_+(-at)+g'_+(at)=2g'_+(at)- \dot y(t)/a,\\
u'(0-,t)&:=&f'_-(-at)+g'_-(at)=2f'_-(-at)+\dot y(t)/a.
\ea
\eeqn
Hence, Eq.~(\ref{4}) writes
$$
m\ddot y(t)=
F(y(t))+2\kappa\,[g'_+(at)-f'_{-}(-at)-\dot y(t)/a],\quad t>0.
$$
Denote
\be\la{1.6}
p(z):=
\frac{u_0(z)+u_0(-z)}{2}+\frac1{2a}\int_{-z}^z u_1(y)\,dy,
\quad z\in\R.
\ee
Therefore, $p(0)=u_0(0)$,
$p'(at)=g'_+(at)-f'_{-}(-at)\in L^2_{{\rm loc}}(\R_+)$, and
we obtain the following evolution equation for $y(t)$, $t>0$:
\beqn
\dot y(t)&=&({a}/{2\kappa})F(y(t))+ a p'(at),\quad t>0,
\quad \mbox{if }\,m=0, \la{1.7'}\\
m\ddot y(t)&=&F(y(t))-({2\kappa}/{a})\dot y(t)+2\kappa p'(at),\quad t>0,
\quad \mbox{if }\,m>0.\la{1.7}
\eeqn
Eq.~(\re{1.1}) implies the following initial condition
for the function $y(t)$:
\be\la{1.8}
y(0)=f_\pm(0)+g_\pm(0)=u_0(0).
\ee
 Eqs (\ref{1.7'}) and (\ref{1.7}) are rewritten in the equivalent integral form,
\beqn\la{1.9}
y(t)&=&\frac a{2\kappa}\int\limits_0^t F(y(s))\,ds+
p(at)-p(0)+y(0),\quad t\geq 0,\quad \mbox{if }\,m=0,\\
 \la{3.7}
my(t)&=&\int\limits_0^t ds\int\limits_0^s
F(y(\tau))\,d\tau+
\frac{2\kappa}{a}\int\limits_0^t (p(as)-y(s))\,ds\nonumber\\
&&+
my(0)+m\dot y(0) t+\frac{2\kappa}{a}\left(y(0)-p(0)\right) t,\quad t\ge0,\quad
\mbox{if }\,m>0.
\eeqn
 Lemma \ref{l2.1} below implies Proposition \ref{p1} immediately.

\begin{lemma}\la{l2.1}
(i) Let $m=0$ and all assumptions of Proposition \ref{p1} hold. Then
for any $y_0\in\R$,
Eq.~(\ref{1.7'}) has a unique solution $y(t)=U(t,0)y_0\in C(\R_+)$.\\
(ii) Let $m>0$. Then
for any $(y_0, y_1)\in\R^2$,
Eq.~(\ref{1.7}) has a unique solution
$(y(t),\dot y(t))=U(t,0)(y_0,y_1)$, and $y(t)\in C^1(\R_+)$.
\\
(iii) For $m\ge0$, the following bound holds,
\be\la{1.10}
\sup_{[0,\tau]}
\left[m|\dot y(t)|+|y(t)|\right]\le C_1\tau+C_2,
\quad \mbox{for any }\,\tau>0.
\ee
\end{lemma}
{\bf Proof}.
We prove Lemma \ref{l2.1} only in the case when $m>0$.
For $m=0$ the proof is similarly.
It follows from (\ref{3.7}), condition (\ref{F}) and
the contraction mapping principle  that for any fixed initial data
$y(0+)$ and $\dot{y}(0+)$, the solution $y(t)$ to  Eq.~(\ref{3.7}) has a unique solution
on a certain interval $t\in [0,\varepsilon)$ with an $\varepsilon$,
$\varepsilon >0$.
Let us derive an a priori estimate for $y(t)$.
This estimate will imply the
existence and uniqueness of the global  solution of (\ref{1.7})
for any $y(0+)$ and $\dot{y}(0+)$.
We multiply Eq.~(\ref{1.7}) by $\dot y(t)$. Using
$\ds\frac{d}{dt}V(y(t))=-F(y(t))\dot y(t)$,
we obtain
$$
\frac{d}{dt}\left(\frac{m\dot y^2(t)}{2}+ V(y(t))\right)=
2\kappa p'(at)\dot y(t)-\frac{2\kappa}{a}\dot y^2(t)
\le \frac{a\kappa}{2} (p'(at))^2.
$$
Let us integrate this inequality and obtain
$$
\frac{m\dot y^2(t)}{2}+ V(y(t))\le \frac{m \dot y^2(0)}{2}+V(y(0))+
\frac{a\kappa}{2} \int_0^t |p'(as)|^2\,ds,\quad t>0.
$$
Hence, for any $\tau>0$, there exist constants $C_1,C_2>0$ such that
\be\label{est}
\sup_{t\in[0,\tau]}\left[
\frac{m\dot y^2(t)}{2}+ V(y(t))\right]\le C_1 \tau+C_2.
\ee
Condition (\ref{V}) implies the estimate (\ref{1.10}).
 Lemma \ref{l2.1} is proved.
\medskip

The following result  follows from
the Gronwall inequality and from a priori estimate (\re{1.10})
(see \cite{KMo}).

\begin{lemma}\la{l2.2}
 Let $m=0$ and $y_1(t)$ and $y_2(t)$ be two solutions of Eq.~(\re{1.7'})
with the initial values
$y_1(0)$ and $y_2(0)$, respectively.
Then for every $\tau>0$,
\be\la{1.14}
 \Vert\dot y_1(t)-\dot y_2(t)\Vert_{L^2(0,\tau)}+
\max_{[0,\tau]} |y_1(t)-y_2(t)|
\le C(\tau)|y_1(0)-y_2(0)|,
\ee
where a constant $C(\tau)$ is bounded for bounded
$y_1(0), y_2(0)$.
The similar result holds for Eq.~(\ref{1.7}) in the case $m\not=0$.
\end{lemma}

\setcounter{equation}{0}
\section{The proof of the main result}\label{sec3}

Since $(u_0,u_1)\in{\cal H}^{\omega}$, the function $p$ defined in (\ref{1.6})
has  the following properties:
$p\in C^1(\R)$,  $p(z\pm\omega)=p(z)$, $\pm z>0$.
Then the function $p'(at)$ in Eqs~(\ref{1.7'}) and (\ref{1.7}) is periodic with
$\omega/a$--period, and $p'(at)\in C(\R_+)$.

\subsection{The string--spring system ($m=0$)}

At first, we study the behavior of solutions to Eq.~(\ref{1.7'}).

\begin{theorem}\la{t2.1}
Let condition (\ref{F}) hold. Then the following assertions are true.\\
(i) All solutions of Eq.~(\ref{1.7'}) are bounded.\\
(ii) Eq.~(\ref{1.7'}) has at least one $\omega/a$-periodic solution. \\
(iii) Any solution $y(t)$ of Eq.~(\re{1.7'}) either is $\omega/a$-periodic
 or tends to an $\omega/a$-periodic solution $y_p(t)$ as $t\to\infty$
such that for every $R>0$,
\be\la{2.4}
\int_t^{t+R} |\dot y(s)-\dot y_p(s)|^2\,ds+
\sup_{s\in[t,t+R]}|y(s)-y_p(s)|\to 0 \quad \mbox{as }\,t\to\infty.
\ee
\end{theorem}

The items (i) and (ii) follows from the results of \cite[\S 9]{P0}.
These assertions imply (iii) by Theorem 9.1 from \cite{P0}.
For completeness of exposition we give the proof of
Theorem \ref{t2.1} in Appendix~A.

{\bf Proof of Theorem \ref{t1}}.
Let $u(x,t)$ be a solution of the problem (\ref{3})--(\ref{4'}), (\ref{2'}).
Then $u(0,t)=y(t)$ is the solution of Eq. (\ref{1.7'}) with the initial condition
$y(0)=u_0(0)$. In Appendix~A we will show that for any $y_0\in\R$
there exists the limit of $U(n\omega/a,0)y_0$ as $n\to\infty$.
Write $\bar y_0:=\lim\limits_{n\to\infty} U(n\omega/a,0)u_0(0)$.
Then $y_p(t)=U(t,0)\bar y_0$ is the $\omega/a$-periodic solution of Eq. (\ref{1.7'})
and convergence (\ref{2.4}) holds (see Appendix~A).

Put $\bar u_0(x)=u_0(x)-u_0(0)+\bar y_0$ and
define functions $\bar f_\pm(x)$ and $\bar g_\pm(x)$
 so as  $f_\pm(x)$ and $g_\pm(x)$ in (\ref{1.1})
but with $\bar u_0(x)$ instead of $u_0(x)$.
Introduce a function $u_p(x,t)$ as follows
\beqn\label{2.11}
u_p(x,t)=\left\{\ba{lll}
\bar f_+(x-at)+\bar g_+(x+at)& \mbox{for }\,\,x\ge at\\
 y_p(t-x/a)+\bar g_+(x+at)-\bar g_+(at-x)& \mbox{for }\,\,0\le x<at\\
 y_p(t+x/a)+\bar f_-(x-at)-\bar f_-(-at-x)&\mbox{for }\,\, -at\le x<0\\
\bar f_-(x-at)+\bar g_-(x+at)& \mbox{for }\,\,x<-at
\ea
\right.
\eeqn
Then $u_p(x,t)$ is the solution of (\ref{3})--(\ref{4'})
with the initial data $(\bar u_0, u_1)$ and $u_p(0,t)=y_p(t)$.
 Since $(\bar u_0,u_1)\in{\cal H}^{\omega}$,
the functions $\bar f_-(\pm x-at)$ and $\bar g_+(\pm x+at)$ in (\ref{2.11})
are $\omega/a$--periodic in $t$.
Then the equality~(\ref{10}) holds,
 and the convergence~(\ref{9}) follows from (\ref{1.11}) and (\ref{2.4}).

\begin{remark}\label{c1.7}
{\rm Let us consider the problem (\ref{1}) for $t>0$ with initial data
(\ref{2'}), satisfying the following conditions:
 $(u_0,u_1)\in{\cal H}$ and $u_1$ has a form
\beqn\label{con_e}
u_1(x)=\left\{
\ba{ll}
a(2p'_+(x)-u'_0(x)),\quad x\ge0,\\
a(u_0'(x)-2p'_-(x)),\quad x<0,
\ea
\right.
\eeqn
where $p_\pm\in C^1(\R_\pm)$ and $p_\pm(x)$ is $\omega$-periodic for $\pm x>0$.
Then $f_-(z)=p_-(z)$ for $z<0$ and $g_+(z)=p_+(z)$ for $z>0$.
Hence, by formula (\ref{1.11}),
the solution $u(x,t)$ for $t>0$ has the form
\beqn\nonumber
u(x,t)=\left\{\ba{lll}
u_0(x-at)-p_+(x-at)+p_+(x+at)& \mbox{for }\,x>at\\
y(t-x/a)+p_+(x+at)-p_+(at-x)& \mbox{for }\,0<x<at\\
y(t+x/a)+p_-(x-at)-p_-(-at-x)&\mbox{for }\, -at<x<0\\
p_-(x-at)+u_0(x+at)-p_-(x+at)& \mbox{for }\,x<-at
\ea
\right.
\eeqn
where $y(t)$ is a solution of Eq.~(\ref{1.7'}) with the
$\omega$-periodic function $p(x):=p_+(x)+p_-(-x)$, $x>0$,
and satisfies the initial condition~(\ref{1.8}).
Then the results of Theorems~\ref{t2.1} and \ref{t1}
hold as $t\to+\infty$.}
\end{remark}

\subsection{The string--oscillator system ($m>0$)}
Put $c=1/m$, $k=2\kappa/(am)=2\sqrt{\kappa\mu}/m$. Then
 Eq.~(\ref{1.7}) is equivalent to the following system
\beqn\label{5.1}
\left\{\ba{rcl}
\dot y&=&v,\\
\dot v&=&c F(y)-k\, v+ka\, p'(at).
\ea
\right.
\eeqn
Denote by
$Y(t,Y_0,t_0)=\left(y(t,Y_0,t_0),\dot y(t,Y_0,t_0)\right)=U(t,t_0)Y_0$
the solution of the Cauchy problem for the system (\ref{5.1})
with the initial data
\be\label{5.2}
Y_0=(y,\dot y)|_{t=t_0}=(y_0,y_1).
\ee
\begin{definition}
The system is called {\it dissipative} (or D-system)
if for any $(Y_0,t_0)\in\R^3$ there exists a $R$, $R>0$, such that
$\lim\limits_{t\to\infty}\Vert Y(t,Y_0,t_0) \Vert<R$.
\end{definition}
\begin{lemma}\label{l3.5}
Let condition (\ref{F}) hold. Then the following assertions hold.\\
(i) The system (\ref{5.1}) is dissipative,
and  there exist constants $M,N>0$
such that for large time the solutions of the system (\ref{5.1})
belong to a bounded set
\be\label{5.3}
\{(y_0,y_1)\in\R^2:\,|y_0|\le M,\,\,|y_1|\le N\},
\ee
and $M$ and $N$ are independent on the parameters $k$ and $c$
of the system (\ref{5.1}).\\
(ii) The system (\ref{5.1}) has at least one $\omega/a$--periodic solution.
\end{lemma}

The item (i) of Lemma \ref{l3.5} follows from the results of Cartwright and Littlewood, Reuter and others
 (see \cite{C, CL, Re} and
the review works \cite[Chapter VII]{SC}, \cite[Chapter XI, \S4]{Le},
 and \cite[Theorem 5.5.4]{RSC}).
According to the Opial theorem (see, e.g., \cite[Theorem 5.3.6]{RSC})
instead of condition~(\ref{F}) it suffices to assume that
$$
\lim_{y\to+\infty} F(y)<-r,\quad
\lim_{y\to-\infty} F(y)>r,\quad \mbox{where }\,r=\max_{t\in\R} |p'(at)|.
$$
Item (i) implies item (ii) by the Brouwer Fixed Point Theorem
  (see \cite[Chapter 1, \S 2]{P0}).
\medskip

Introduce a mapping $T:\R^2\to\R^2$  as $T=U(\omega_0,0)$, $\omega_0:=\omega/a$.
The map $T$ is called {\it the Poincar\'e transformation}
associated with the periodic system (\ref{5.1}).
Lemma \ref{l3.5} and the Pliss results (see \cite[Chapter 2, \S2]{P}) imply that
there exists an invariant set $I$ w.r.t. $T$, i.e., $TI=I$.
This set is called {\it characteristical set} of the dissipative system
(\ref{5.1}) or a {\it global attractor} of the diffeomorphism $T$.
The set $I$ has the following properties  (see \cite{P0}--\cite{RSC}):
\begin{itemize}
\item
 $I$ is closed and bounded.
\item
$I$ is {\it stable w.r.t.} $T$, i.e., for any $\varepsilon>0$
there exists $\delta>0$ such that if $\rho(Y_0,I)<\delta$ then
$\rho(T^mY_0,I)<\varepsilon$ for every $m\in\N$.
\item
For all $Y_0\in\R^2$, $\rho(T^nY_0,I)\to0$, $n\to\infty$.
\item
There exists a fixed point of the mapping $T$ belonging to $I$, i.e.,
there exists an $\omega_0=\omega/a$-periodic solution (or {\it harmonics}) of the system (\ref{5.1}).
\item
The set $I$ has zero Lebesgue measure by Theorem 1.9 from \cite{P}.
\end{itemize}

Define a set ${\cal S}$ as
$$
{\cal S}:=\{(Y,t)\in\R^3:\, Y=Y(t,Y_0,t_0),\,\, Y_0\in I,\,\,t\in\R\}.
$$
The set ${\cal S}$ has the following properties:
\begin{itemize}
\item
 ${\cal S}$ is bounded and closed.
\item
 ${\cal S}$ is $\omega_0$--periodic, i.e.,
for $(Y,t)\in {\cal S}$, $(Y,t+n\omega_0)\in {\cal S}$, $\forall n\in\N$.
\item
${\cal S}$ is invariant, i.e., 
if $(Y_0,t_0)\in {\cal S}$, then
$(Y(t,Y_0,t_0),t)\in {\cal S}$ for all $t\ge t_0$.
\item
 ${\cal S}$ is stable, i.e., $\forall\ve>0$
$\exists \delta>0$ such that if
$\rho(Y_0,{\cal S}_{t_0})<\delta$, then
$\rho(Y(t,Y_0,t_0),{\cal S}_{t})<\ve$, $\forall t\ge t_0$,
where ${\cal S}_{\tau}={\cal S}\cap\{t=\tau\}$.
\item
 ${\cal S}$ is stable in whole, i.e.,
for all $Y(t,Y_0,t_0)\in\R^2$ we have $\lim\limits_{t\to\infty}
\rho(Y(t,Y_0,t_0),{\cal S}_{t})=0$.
\end{itemize}

However, these properties of ${\cal S}$ do not imply, in general,
 the convergence (\ref{2.4}).
Now we consider the particular case of the system (\ref{5.1})
when $I$ has a unique point.
Then (\ref{5.1}) is called
{\it the system with convergence} (see \cite[\S 7, Definition 7.1]{P0}).
In this case,
the system (\ref{5.1})
has a unique stable $\omega_0$-periodic solution $Y_p(t)$,
and any another solution $Y(t,Y_0,t_0)$ tends to this periodic solution, i.e.,
$\lim\limits_{t\to\infty}\Vert Y(t,Y_0,t_0)-Y_p(t)\Vert=0$, and
the result (\ref{2.4}) follows.
\medskip

Below we give examples of the restrictions on the function
$F(y)$ when the system (\ref{5.1}) has convergence property.
\begin{example}\la{ex1}
{\rm Assume that
\begin{description}
  \item[(F1)] $F(y)=-ry$ with a constant $r>0$.
\end{description}
Then by the Levinson theorem (see \cite{Lev}, \cite[Theorem 8.1]{P0}, \cite[Theorem~5.2.1]{RSC}),
 Eq.~(\ref{1.7}) has a unique $\omega_0$--periodic
solution and all other solutions tend to this periodic solution as $t\to+\infty$.}
\end{example}
\begin{example}\la{ex2}
{\rm Assume that for $y_1\not=y_2$, we have
\begin{description}
  \item[(F2)]
$ k^2/2-1\le\ds -c\frac{F(y_2)-F(y_1)}{y_2-y_1}\le 1,\quad 1<k^2/2\le 2$,\\
 where $k=\sqrt{\kappa\mu}/m$ is the constant in (\ref{5.1}).
\end{description}
Then according to the Zlam\'al theorem (see, e.g., \cite[Theorem 5.3.2]{RSC})
 all solutions tend  exponentially
to a unique periodic solution as $t\to+\infty$.}
\end{example}
\begin{example}\la{ex3}
{\rm (see \cite[Theorem 8.4]{P0}, \cite[Ch.XI, \S 5]{Le} or \cite{Sh})
Assume that
\begin{description}
  \item[(F3)]
 $F\in C^2(\R)$,$F'(y)<0$ for $|y|\le M$;
 $\exists \beta>0$ such that $F(y)\,{\rm sgn}\, y\le-\beta$ for $|y|\ge M$,
with the constant $M$ from the bound~(\ref{5.3}).
Moreover, the constant $k$ from (\ref{5.1}) is enough large,
$$
k>(1/2) N\max_{|y|\le M}\left({|F''(y)|}/{|F'(y)|}\right),
$$
where the constant $N$ is defined in the bound~(\ref{5.3}).
\end{description}
Then the system (\ref{5.1}) has convergence property.
For instance, the function $F(y)=-ay^3-by$ with constants $a,b>0$
satisfies these conditions.}
\end{example}

Note that condition {\bf(F1)} is a particular case of {\bf(F3)}.
\begin{cor}
Let condition {\bf(F2)} or {\bf(F3)}   be true.
Then the following assertions hold.\\
(i) There exists a unique $\omega_0$-periodic solution $y_p(t)$ of Eq. (\ref{1.7}),
and for any another solution $y(t)$ the convergence (\ref{2.4}) holds.\\
(ii) The convergence (\ref{9}) holds with the function $u_p(x,t)$ satisfying (\ref{10}).
\end{cor}

The assertion (i) follows from the results mentioned above.
Now we check item (ii). Indeed,
let $u(x,t)$ be a solution of the problem (\ref{3})--(\ref{4'}), (\ref{2'}), (\ref{2''}).
Then there exists
$\lim\limits_{n\to\infty} T^n(u_0(0),y_1)=:(\bar y_0,\bar y_1)$
and $(\bar y_0,\bar y_1)$ is a unique point of the set $I$.
Hence $(y_p(t),\dot y_p(t)) =U(t,0)(\bar y_0,\bar y_1)$ is the unique $\omega/a$-periodic solution
of the system (\ref{5.1}) and convergence (\ref{2.4}) holds.
Put $\bar u_0(x)=u_0(x)-u_0(0)+\bar y_0$ and
define functions $\bar f_\pm(x)$ and $\bar g_\pm(x)$  by formulas (\ref{1.1})
but with $\bar u_0(x)$ instead of $u_0(x)$.
Define $u_p(x,t)$  by (\ref{2.11}).
Then $u_p(x,t)$ is the solution of the problem (\ref{3})--(\ref{4'})
 with the initial data $(\bar u_0, u_1, \bar y_1)$.
 Since $(\bar u_0,u_1)\in{\cal H}^{\omega}$,
the functions $\bar f_-(\pm x-at)$ and $\bar g_+(\pm x+at)$ in (\ref{2.11})
are $\omega/a$--periodic in $t$.
Hence 
 the equality~(\ref{10}) holds,
 and the convergence~(\ref{9}) follows from (\ref{1.11}) and (\ref{2.4}).

\setcounter{equation}{0}
\section{Appendix A: Proof of Theorem 3.1}
For simplicity, instead of Eq.~(\ref{1.7'})
we consider the following equation:
\be\label{4.0}
\dot y(t)=F(y(t))+P(t),\quad t>0,
\ee
where $P(t)\in C(0,+\infty)$ is a periodic function with period
$\omega_0=\omega/a$,
the function $F(y)$ satisfies the condition (\ref{F}).
Write
\be\la{q}
q=\max_{t\in\R}|P(t)|<\infty.
\ee
Denote by $U(t,s)$, $t\ge s$,
a solving operator of the following Cauchy problem
for the function $y(t)$:
\beqn
\dot y(t)&=&F(y(t))+P(t),\quad t>s,    \la{2.2}\\
y|_{t=s}&=&y_0.                    \la{2.3}
\eeqn
Then
$U(t,s):\R^1\to\R^1$ transforms the initial condition $y_0\in\R^1$ for $t=s$
to the solution  $y(t)$ of the problem (\ref{2.2})--(\ref{2.3}) in time $t$:
$U(t,s):y_0\to y(t)\equiv y(t,y_0,s)$.
\begin{remark}\label{r2.1}
(i) $U(t,s) U(s,r)=U(t,r)$ for $t>s>r$,
(ii) $U(s,s)=I_d$, where $I_d$ is identity operator on $\R$,
(iii) $U(t+\omega_0,s+\omega_0)=U(t,s)$, $t,s\in\R$.
\end{remark}

At first, we prove the following lemma.
\begin{lemma}\la{l2.21}
 There exists a bounded interval $B=[y_-,y_+]\subset\R^1$
such that $U(t,s) :B\subset B$, $t\ge s$.
Moreover, all solutions $y(t)$,  $t\ge s$,
 to problem (\ref{2.2})--(\ref{2.3}) for finite time come in the region $B$.
\end{lemma}
{\bf Proof}.
It follows from (\ref{F}) and (\ref{q})
that there exist points $y_-<0$ and $y_+>0$
such that
\beqn
F(y)-q>0 &\quad \mbox{for }& y\le y_-;\nonumber\\
F(y)+q<0 &\quad \mbox{for }& y\ge y_+.\nonumber
\eeqn
Then by Eq.~(\ref{2.2}) we have
$F(y(t))-q\le \dot y(t)\le F(y(t))+q$.
Therefore, if $y(t)\ge y_+$, then $\dot y(t)<0$,
if  $y(t)\le y_-$, then  $\dot y(t)>0$.
Lemma \ref{l2.21} is proved.
\medskip

Denote by $T:=U(\omega_0,0)$ the Poincar\'e transformation
associated with Eq.~(\ref{2.2}).
Since  $F\in C^1(\R)$, the right hand side of
Eq.~(\ref{2.2}) is continuously differentiable on $y$.
Hence $T$ is  continuously differentiable mapping
$\R$ on $\R$.
It is easy to verify that
there exists a continuously differentiable mapping
$T^{-1}=U(0,\omega_0)$, i.e., $T$ is a diffeomorphism $\R$ on $\R$.
It follows from Lemma \ref{l2.21}
 and the Brouwer Fixed Point Theorem that the
mapping $T$ has a fixed point belonging to the interval $B$
(see \cite[Theorem 2.9.1]{RSC} or Massera's theorem \cite{Mas}).

Denote by ${\cal Z}$ a set of fixed points of the mapping $T$,
$$
{\cal Z}=\{y\in\R:\,Ty=y\}.
$$
 \begin{lemma}\la{l2.22}
 For every  $y_0\in\R$, there exists a $\bar y_0\in{\cal Z}$
such that $T^n y_0\to \bar y_0$ as $n\to\infty$.
\end{lemma}
{\bf Proof}. Since $T$ is a diffeomorphism, then $T$ is a monotone function.
If $y>y_+$, then $Ty<y$, and if $y<y_-$ then $Ty>y$.
Let $y_0$ be an arbitrary point.
It follows from Lemma~\ref{l2.21}
that $\exists N(y_0)\in\N$ such that
$\forall n>N(y_0)$ we have $T^n y_0\in B$.
Further, we assume that $y_0\in B$.
It is possible 3 cases:
\beqn
&(i)& y_0\in (y_-,z_-),\quad \mbox{where }~~z_-=\inf_{z\in {\cal Z}}z,
\nonumber\\
&(ii)& y_0\in (z_+,y_+),\quad \mbox{where }~~z_+=\sup_{z\in {\cal Z}}z,
\nonumber\\
&(iii)& y_0\in (z_i,z_{i+1}),\quad \mbox{where }~~z_i,~z_{i+1}\in
{\cal Z}~~\mbox{are two neighboring fixed points of }\,T.
\nonumber
\eeqn
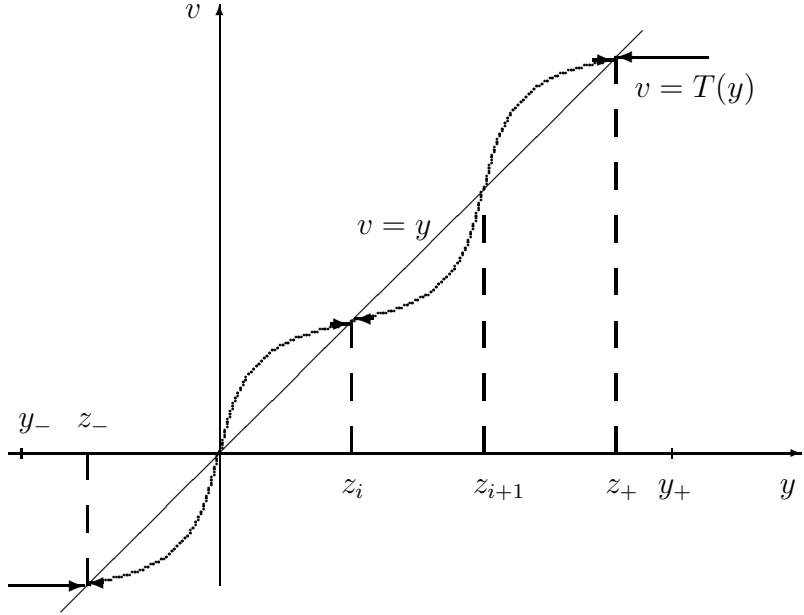
\begin{picture}(500,300)
\put(100,50){\begin{picture}(500,300)
\put(50,0){\vector(0,1){220}}
\put(-30,50){\vector(1,0){300}}
\put(-10,-10){\line(1,1){220}}
\put(  0,  0){\circle*{1}}
\put(  1,  0){\circle*{1}}
\put(  2,  1){\circle*{1}}
\put(  3,  1){\circle*{1}}
\put(  4,  1){\circle*{1}}
\put(  4,  1){\circle*{1}}
\put(  5,  1){\circle*{1}}
\put(  6,  2){\circle*{1}}
\put(  7,  2){\circle*{1}}
\put(  8,  2){\circle*{1}}
\put(  8,  2){\circle*{1}}
\put(  9,  2){\circle*{1}}
\put( 10,  3){\circle*{1}}
\put( 11,  3){\circle*{1}}
\put( 12,  3){\circle*{1}}
\put( 12,  3){\circle*{1}}
\put( 13,  3){\circle*{1}}
\put( 14,  4){\circle*{1}}
\put( 15,  4){\circle*{1}}
\put( 15,  4){\circle*{1}}
\put( 16,  4){\circle*{1}}
\put( 17,  5){\circle*{1}}
\put( 18,  5){\circle*{1}}
\put( 18,  5){\circle*{1}}
\put( 19,  5){\circle*{1}}
\put( 20,  6){\circle*{1}}
\put( 21,  6){\circle*{1}}
\put( 21,  6){\circle*{1}}
\put( 22,  7){\circle*{1}}
\put( 23,  7){\circle*{1}}
\put( 23,  7){\circle*{1}}
\put( 24,  8){\circle*{1}}
\put( 25,  8){\circle*{1}}
\put( 25,  8){\circle*{1}}
\put( 26,  9){\circle*{1}}
\put( 27,  9){\circle*{1}}
\put( 27,  9){\circle*{1}}
\put( 28, 10){\circle*{1}}
\put( 28, 10){\circle*{1}}
\put( 29, 10){\circle*{1}}
\put( 30, 11){\circle*{1}}
\put( 30, 11){\circle*{1}}
\put( 31, 12){\circle*{1}}
\put( 31, 12){\circle*{1}}
\put( 32, 13){\circle*{1}}
\put( 32, 13){\circle*{1}}
\put( 33, 14){\circle*{1}}
\put( 33, 14){\circle*{1}}
\put( 34, 15){\circle*{1}}
\put( 34, 15){\circle*{1}}
\put( 35, 16){\circle*{1}}
\put( 35, 16){\circle*{1}}
\put( 36, 17){\circle*{1}}
\put( 36, 17){\circle*{1}}
\put( 37, 18){\circle*{1}}
\put( 37, 18){\circle*{1}}
\put( 38, 19){\circle*{1}}
\put( 38, 19){\circle*{1}}
\put( 39, 20){\circle*{1}}
\put( 39, 20){\circle*{1}}
\put( 39, 21){\circle*{1}}
\put( 40, 22){\circle*{1}}
\put( 40, 22){\circle*{1}}
\put( 41, 23){\circle*{1}}
\put( 41, 24){\circle*{1}}
\put( 41, 24){\circle*{1}}
\put( 42, 25){\circle*{1}}
\put( 42, 26){\circle*{1}}
\put( 42, 26){\circle*{1}}
\put( 43, 27){\circle*{1}}
\put( 43, 28){\circle*{1}}
\put( 43, 28){\circle*{1}}
\put( 44, 29){\circle*{1}}
\put( 44, 30){\circle*{1}}
\put( 44, 30){\circle*{1}}
\put( 44, 31){\circle*{1}}
\put( 45, 32){\circle*{1}}
\put( 45, 33){\circle*{1}}
\put( 45, 33){\circle*{1}}
\put( 45, 34){\circle*{1}}
\put( 46, 35){\circle*{1}}
\put( 46, 36){\circle*{1}}
\put( 46, 36){\circle*{1}}
\put( 46, 37){\circle*{1}}
\put( 47, 38){\circle*{1}}
\put( 47, 39){\circle*{1}}
\put( 47, 40){\circle*{1}}
\put( 47, 40){\circle*{1}}
\put( 47, 41){\circle*{1}}
\put( 48, 42){\circle*{1}}
\put( 48, 43){\circle*{1}}
\put( 48, 44){\circle*{1}}
\put( 48, 44){\circle*{1}}
\put( 48, 45){\circle*{1}}
\put( 49, 46){\circle*{1}}
\put( 49, 47){\circle*{1}}
\put( 49, 48){\circle*{1}}
\put( 49, 48){\circle*{1}}
\put( 49, 49){\circle*{1}}
\put( 50, 50){\circle*{1}}
\put( 50, 50){\circle*{1}}
\put( 50, 51){\circle*{1}}
\put( 51, 52){\circle*{1}}
\put( 51, 53){\circle*{1}}
\put( 51, 54){\circle*{1}}
\put( 51, 54){\circle*{1}}
\put( 51, 55){\circle*{1}}
\put( 52, 56){\circle*{1}}
\put( 52, 57){\circle*{1}}
\put( 52, 58){\circle*{1}}
\put( 52, 58){\circle*{1}}
\put( 52, 59){\circle*{1}}
\put( 53, 60){\circle*{1}}
\put( 53, 61){\circle*{1}}
\put( 53, 62){\circle*{1}}
\put( 53, 62){\circle*{1}}
\put( 53, 63){\circle*{1}}
\put( 54, 64){\circle*{1}}
\put( 54, 65){\circle*{1}}
\put( 54, 65){\circle*{1}}
\put( 54, 66){\circle*{1}}
\put( 55, 67){\circle*{1}}
\put( 55, 68){\circle*{1}}
\put( 55, 68){\circle*{1}}
\put( 55, 69){\circle*{1}}
\put( 56, 70){\circle*{1}}
\put( 56, 71){\circle*{1}}
\put( 56, 71){\circle*{1}}
\put( 57, 72){\circle*{1}}
\put( 57, 73){\circle*{1}}
\put( 57, 73){\circle*{1}}
\put( 58, 74){\circle*{1}}
\put( 58, 75){\circle*{1}}
\put( 58, 75){\circle*{1}}
\put( 59, 76){\circle*{1}}
\put( 59, 77){\circle*{1}}
\put( 59, 77){\circle*{1}}
\put( 60, 78){\circle*{1}}
\put( 60, 78){\circle*{1}}
\put( 60, 79){\circle*{1}}
\put( 61, 80){\circle*{1}}
\put( 61, 80){\circle*{1}}
\put( 62, 81){\circle*{1}}
\put( 62, 81){\circle*{1}}
\put( 63, 82){\circle*{1}}
\put( 63, 82){\circle*{1}}
\put( 64, 83){\circle*{1}}
\put( 64, 83){\circle*{1}}
\put( 65, 84){\circle*{1}}
\put( 65, 84){\circle*{1}}
\put( 66, 85){\circle*{1}}
\put( 66, 85){\circle*{1}}
\put( 67, 86){\circle*{1}}
\put( 67, 86){\circle*{1}}
\put( 68, 87){\circle*{1}}
\put( 68, 87){\circle*{1}}
\put( 69, 88){\circle*{1}}
\put( 69, 88){\circle*{1}}
\put( 70, 89){\circle*{1}}
\put( 70, 89){\circle*{1}}
\put( 71, 89){\circle*{1}}
\put( 72, 90){\circle*{1}}
\put( 72, 90){\circle*{1}}
\put( 73, 91){\circle*{1}}
\put( 74, 91){\circle*{1}}
\put( 74, 91){\circle*{1}}
\put( 75, 92){\circle*{1}}
\put( 76, 92){\circle*{1}}
\put( 76, 92){\circle*{1}}
\put( 77, 93){\circle*{1}}
\put( 78, 93){\circle*{1}}
\put( 78, 93){\circle*{1}}
\put( 79, 94){\circle*{1}}
\put( 80, 94){\circle*{1}}
\put( 80, 94){\circle*{1}}
\put( 81, 94){\circle*{1}}
\put( 82, 95){\circle*{1}}
\put( 83, 95){\circle*{1}}
\put( 83, 95){\circle*{1}}
\put( 84, 95){\circle*{1}}
\put( 85, 96){\circle*{1}}
\put( 86, 96){\circle*{1}}
\put( 86, 96){\circle*{1}}
\put( 87, 96){\circle*{1}}
\put( 88, 97){\circle*{1}}
\put( 89, 97){\circle*{1}}
\put( 90, 97){\circle*{1}}
\put( 90, 97){\circle*{1}}
\put( 91, 97){\circle*{1}}
\put( 92, 98){\circle*{1}}
\put( 93, 98){\circle*{1}}
\put( 94, 98){\circle*{1}}
\put( 94, 98){\circle*{1}}
\put( 95, 98){\circle*{1}}
\put( 96, 99){\circle*{1}}
\put( 97, 99){\circle*{1}}
\put( 98, 99){\circle*{1}}
\put( 98, 99){\circle*{1}}
\put( 99, 99){\circle*{1}}
\put(100,100){\circle*{1}}
\put(100,100){\circle*{1}}
\put(101,100){\circle*{1}}
\put(102,101){\circle*{1}}
\put(103,101){\circle*{1}}
\put(104,101){\circle*{1}}
\put(104,101){\circle*{1}}
\put(105,101){\circle*{1}}
\put(106,102){\circle*{1}}
\put(107,102){\circle*{1}}
\put(108,102){\circle*{1}}
\put(108,102){\circle*{1}}
\put(109,102){\circle*{1}}
\put(110,103){\circle*{1}}
\put(111,103){\circle*{1}}
\put(112,103){\circle*{1}}
\put(112,103){\circle*{1}}
\put(113,103){\circle*{1}}
\put(114,104){\circle*{1}}
\put(115,104){\circle*{1}}
\put(115,104){\circle*{1}}
\put(116,104){\circle*{1}}
\put(117,105){\circle*{1}}
\put(118,105){\circle*{1}}
\put(118,105){\circle*{1}}
\put(119,105){\circle*{1}}
\put(120,106){\circle*{1}}
\put(121,106){\circle*{1}}
\put(121,106){\circle*{1}}
\put(122,107){\circle*{1}}
\put(123,107){\circle*{1}}
\put(123,107){\circle*{1}}
\put(124,108){\circle*{1}}
\put(125,108){\circle*{1}}
\put(125,108){\circle*{1}}
\put(126,109){\circle*{1}}
\put(127,109){\circle*{1}}
\put(127,109){\circle*{1}}
\put(128,110){\circle*{1}}
\put(128,110){\circle*{1}}
\put(129,110){\circle*{1}}
\put(130,111){\circle*{1}}
\put(130,111){\circle*{1}}
\put(131,112){\circle*{1}}
\put(131,112){\circle*{1}}
\put(132,113){\circle*{1}}
\put(132,113){\circle*{1}}
\put(133,114){\circle*{1}}
\put(133,114){\circle*{1}}
\put(134,115){\circle*{1}}
\put(134,115){\circle*{1}}
\put(135,116){\circle*{1}}
\put(135,116){\circle*{1}}
\put(136,117){\circle*{1}}
\put(136,117){\circle*{1}}
\put(137,118){\circle*{1}}
\put(137,118){\circle*{1}}
\put(138,119){\circle*{1}}
\put(138,119){\circle*{1}}
\put(139,120){\circle*{1}}
\put(139,120){\circle*{1}}
\put(139,121){\circle*{1}}
\put(140,122){\circle*{1}}
\put(140,122){\circle*{1}}
\put(141,123){\circle*{1}}
\put(141,124){\circle*{1}}
\put(141,124){\circle*{1}}
\put(142,125){\circle*{1}}
\put(142,126){\circle*{1}}
\put(142,126){\circle*{1}}
\put(143,127){\circle*{1}}
\put(143,128){\circle*{1}}
\put(143,128){\circle*{1}}
\put(144,129){\circle*{1}}
\put(144,130){\circle*{1}}
\put(144,130){\circle*{1}}
\put(144,131){\circle*{1}}
\put(145,132){\circle*{1}}
\put(145,133){\circle*{1}}
\put(145,133){\circle*{1}}
\put(145,134){\circle*{1}}
\put(146,135){\circle*{1}}
\put(146,136){\circle*{1}}
\put(146,136){\circle*{1}}
\put(146,137){\circle*{1}}
\put(147,138){\circle*{1}}
\put(147,139){\circle*{1}}
\put(147,140){\circle*{1}}
\put(147,140){\circle*{1}}
\put(147,141){\circle*{1}}
\put(148,142){\circle*{1}}
\put(148,143){\circle*{1}}
\put(148,144){\circle*{1}}
\put(148,144){\circle*{1}}
\put(148,145){\circle*{1}}
\put(149,146){\circle*{1}}
\put(149,147){\circle*{1}}
\put(149,148){\circle*{1}}
\put(149,148){\circle*{1}}
\put(149,149){\circle*{1}}
\put(150,150){\circle*{1}}
\put(150,150){\circle*{1}}
\put(150,151){\circle*{1}}
\put(151,152){\circle*{1}}
\put(151,153){\circle*{1}}
\put(151,154){\circle*{1}}
\put(151,154){\circle*{1}}
\put(151,155){\circle*{1}}
\put(152,156){\circle*{1}}
\put(152,157){\circle*{1}}
\put(152,158){\circle*{1}}
\put(152,158){\circle*{1}}
\put(152,159){\circle*{1}}
\put(153,160){\circle*{1}}
\put(153,161){\circle*{1}}
\put(153,162){\circle*{1}}
\put(153,162){\circle*{1}}
\put(153,163){\circle*{1}}
\put(154,164){\circle*{1}}
\put(154,165){\circle*{1}}
\put(154,165){\circle*{1}}
\put(154,166){\circle*{1}}
\put(155,167){\circle*{1}}
\put(155,168){\circle*{1}}
\put(155,168){\circle*{1}}
\put(155,169){\circle*{1}}
\put(156,170){\circle*{1}}
\put(156,171){\circle*{1}}
\put(156,171){\circle*{1}}
\put(157,172){\circle*{1}}
\put(157,173){\circle*{1}}
\put(157,173){\circle*{1}}
\put(158,174){\circle*{1}}
\put(158,175){\circle*{1}}
\put(158,175){\circle*{1}}
\put(159,176){\circle*{1}}
\put(159,177){\circle*{1}}
\put(159,177){\circle*{1}}
\put(160,178){\circle*{1}}
\put(160,178){\circle*{1}}
\put(160,179){\circle*{1}}
\put(161,180){\circle*{1}}
\put(161,180){\circle*{1}}
\put(162,181){\circle*{1}}
\put(162,181){\circle*{1}}
\put(163,182){\circle*{1}}
\put(163,182){\circle*{1}}
\put(164,183){\circle*{1}}
\put(164,183){\circle*{1}}
\put(165,184){\circle*{1}}
\put(165,184){\circle*{1}}
\put(166,185){\circle*{1}}
\put(166,185){\circle*{1}}
\put(167,186){\circle*{1}}
\put(167,186){\circle*{1}}
\put(168,187){\circle*{1}}
\put(168,187){\circle*{1}}
\put(169,188){\circle*{1}}
\put(169,188){\circle*{1}}
\put(170,189){\circle*{1}}
\put(170,189){\circle*{1}}
\put(171,189){\circle*{1}}
\put(172,190){\circle*{1}}
\put(172,190){\circle*{1}}
\put(173,191){\circle*{1}}
\put(174,191){\circle*{1}}
\put(174,191){\circle*{1}}
\put(175,192){\circle*{1}}
\put(176,192){\circle*{1}}
\put(176,192){\circle*{1}}
\put(177,193){\circle*{1}}
\put(178,193){\circle*{1}}
\put(178,193){\circle*{1}}
\put(179,194){\circle*{1}}
\put(180,194){\circle*{1}}
\put(180,194){\circle*{1}}
\put(181,194){\circle*{1}}
\put(182,195){\circle*{1}}
\put(183,195){\circle*{1}}
\put(183,195){\circle*{1}}
\put(184,195){\circle*{1}}
\put(185,196){\circle*{1}}
\put(186,196){\circle*{1}}
\put(186,196){\circle*{1}}
\put(187,196){\circle*{1}}
\put(188,197){\circle*{1}}
\put(189,197){\circle*{1}}
\put(190,197){\circle*{1}}
\put(190,197){\circle*{1}}
\put(191,197){\circle*{1}}
\put(192,198){\circle*{1}}
\put(193,198){\circle*{1}}
\put(194,198){\circle*{1}}
\put(194,198){\circle*{1}}
\put(195,198){\circle*{1}}
\put(196,199){\circle*{1}}
\put(197,199){\circle*{1}}
\put(198,199){\circle*{1}}
\put(198,199){\circle*{1}}
\put(199,199){\circle*{1}}
\put(200,200){\circle*{1}}
\put(94,35){$z_i$}
\put(-6,60){$z_-$}
\put(144,35){$z_{i+1}$}
\put(-28,60){$y_-$}
\put(194,35){$z_+$}
\put(214,35){$y_+$}
\put(-25,48){\line(0,1){4}}
\put(221,48){\line(0,1){4}}
\put(260,35){$y$}
\put(100,133){$v=y$}
\put(35,215){$v$}
\thicklines
\put(92,99){\vector(1,0){7}}
\put(108,101){\vector(-1,0){7}}
\put(6,1){\vector(-1,0){6}}
\put(235,200){\vector(-1,0){35}}
\put(205,185){$v=T(y)$}
\put(-30,0){\vector(1,0){29}}
\put(191,199){\vector(1,0){8}}
\multiput(100,50)(0,20){3}{\line(0,1){10}}
\multiput(150,50)(0,20){5}{\line(0,1){10}}
\multiput(200,50)(0,20){8}{\line(0,1){10}}
\multiput(0,0)(0,20){3}{\line(0,1){10}}
\end{picture}}
\end{picture}

 In case (i), since $z_-$ is extreme left fixed point of $T$, $Ty>y$
for all $y\in (y_-,z_-)$.
Write $y_n=T^n y_0$. Obviously,  $y_n>y_{n-1}>...>y_0$,
i.e., $y_n$ is an increasing sequence bounded above by $z_-$.
Hence, the limit holds,
$$
\lim_{n\to\infty}y_n=\sup_{n\in\N}y_n=y_s,\quad y_s\in (y_-,z_-],
$$
i.e., $T^n y_0\uparrow y_s$, $n\to\infty$.
Hence, $T^{n+1} y_0\to T y_s=y_s$, $n\to\infty$.
Therefore, $y_s$ is a fixed point of $T$.
Since $z_-$ is extreme left fixed point
of $T$, then $y_s=z_-$.

In case (ii), $Ty<y$, $\forall y\in (z_+,y_+)$.
Hence $y_n=T^ny_0$ is a decreasing
sequence bounded below by $z_+$.
Moreover, $y_n\downarrow z_+$ as $n\to\infty$.

In the case (iii), we have either
(a) $Ty>y$, $y\in (z_i,z_{i+1})$ or
(b) $Ty<y$, $y\in (z_i,z_{i+1})$.
In the case (a) (see the case (i))
$T^ny_0\uparrow z_{i+1}$ as $n\to\infty$.
In the case (b) (see the case (ii))
$T^ny_0\downarrow z_{i}$ as $n\to\infty$.
Theorem \ref{l2.21} is proved.
\medskip\\
{\bf Proof of Theorem \ref{t2.1}}.
It follows from Remark \ref{r2.1} that
$U(n\omega_0,0)=(U(\omega_0,0))^n=T^n$.
Theorem \ref{l2.21} implies that
$\forall y_0\in\R$,
$U(n\omega_0,0)y_0\to \bar y_0$ as $n\to\infty$,
where $\bar y_0$ is a fixed point of $T$.

For $t>0$ we choose
$n\in\N$ such that 
 $t=n\omega_0+\tau$,
where $\tau\in [0,\omega_0)$.
Since $\bar y_0$  is a fixed point of $T$, we have
$$
| U(t,0)y_0-U(t,0)\bar y_0| =|U(\tau,0)T^n y_0-U(\tau,0)T^n \bar y_0|=
| U(\tau,0)T^n y_0-U(\tau,0)\bar y_0|.
$$
Hence,
it follows from Lemma \ref{l2.2} and Lemma \ref{l2.22}
that there exists a constant $C<\infty$
such that for enough large value $n\in \N$, we have
$$
\sup_{\tau\in [0,\omega_0]}
|U(\tau,0)T^ny_0-U(\tau,0)\bar y_0|\le C|T^ny_0-\bar y_0| <\varepsilon.
$$
Denote by $y(t)=U(t,0)y_0$ and $y_p(t)=U(t,0)\bar y_0$ the solutions of the problem
(\ref{2.2})--(\ref{2.3}) with the initial data $y_0$ and $\bar y_0$, respectively.
Note that $y_p(t)$ is a periodic solution of Eq.~(\ref{2.2}),
 because $\bar y_0$ is a fixed point of $T$.
Therefore, for any solution $y(t)$ of Eq.~(\ref{2.2})
 we have $|y(t)-y_p(t)|<\varepsilon$ for enough large value $t>0$.
The convergence~(\ref{2.4}) follows
from condition~(\ref{F}) and Eq.~(\ref{2.2}).

\setcounter{equation}{0}
\section{Appendix B: Limit amplitude principle}
Here we apply the results to the following problem
 for a function $u(x,t)\in C(\R^{2})$:
\beqn
(\mu+m\delta(x))\ddot u(x,t)&=&
\kappa u''(x,t)+\delta(x)F(u(x,t)),\quad t>0,\quad x\in \R,\label{A.1}\\
u(x,t)\mid_{t\le0}&=&p(x+at),\quad x\in \R. \label{A.2}
\eeqn
Here  $m\ge0$,  $a=\sqrt{\kappa/\mu}$.
In the case $m>0$ we assume that condition {\bf (F2)} or {\bf (F3)} holds.
The function $p$ from Eq.~(\ref{A.2}) satisfies the following conditions:
\begin{description}
   \item[P1] $p\in C^1(\R)$.
    \item[P2] There exist numbers $\omega>0$ and $p_0\in\R$ such that $F(p_0)=0$ and
$$
p(z+\omega)=p(z)\,\,\, \mbox{for } z>0,\quad
p(z)=p_0\,\,\,\mbox{for } z\leq 0.
$$
 \end{description}

Note that the function  $p(x+at)$ is a solution
of Eq.~(\re{A.1}) for $t<0$. Therefore, we can consider
Eq.~(\ref{A.1})  for $t\in\R$.
In particular, we have
\beqn\nonumber
u_0(x)=u|_{t=0}=p(x),\quad u_1(x)=\dot u|_{t=0}=ap'(x),\quad x\in\R,
\quad y(0)=u_0(0)=p_0,\quad \dot y(0)=0.
\eeqn
Then $f_\pm(z)=0$ and $g_\pm(z)=p(z)$ for $\pm z>0$.
Therefore, by (\ref{1.11}), 
\beqn\la{1.11'}
u(x,t)=\left\{\begin{array}{lcl}
p(x+at)&  \mbox{for }& \quad x>at,\\
y(t-x/a)-p(at-x)+p(x+at) &  \mbox{for }& \quad 0<x<at,\\
y(t+x/a)  & \mbox{for }&-at<x<0,\\
p_0& \mbox{for }&x<-at.
\end{array}\right.
\eeqn
where $y(t)$ is a solution to Eq.~(\ref{1.7}) (or Eq. (\ref{1.7'})) for $t>0$,
and $y(t)=p_0$ for $t\le0$.
By Proposition~\ref{p1}, the Cauchy problem (\ref{A.1})--(\ref{A.2})
has a unique solution $u(x,t)\in {\cal E}$
for every function $p\in C^1(\R)$.

Let $m\ge0$ and $y_p(t)$ be the $\omega/a$-periodic solution of Eq. (\ref{1.7'})
or (\ref{1.7})
with the initial date $\bar p_0=\lim_{n\to\infty} T^n p_0$ (if $m=0$) or
with the initial data $(\bar p_0,\bar y_1)=\lim_{n\to\infty} T^n (p_0,0)$ (if $m>0$).
We extend $y_p(t)\equiv \bar p_0$ for $t<0$ and define
\beqn\la{2.5}
u_p(x,t)=\left\{\begin{array}{lcl}
y_p(t-x/a)-p(at-x)+p(x+at) & \mbox{for }&x>0,\quad t>0,\\
y_p(t+x/a)    &  \mbox{for }&x<0,\quad t>0.
\end{array}\right.
\eeqn
Then $u_p(x,t)\in{\cal E}$, $u_p(x,t)$ is the solution of Eq. (\re{A.1})
under the condition $u_p(x,t)|_{t\le0}=\bar p(x+at)$,
where $\bar p(x)=\bar p_0+p(x)-p_0$ for $x\in\R$.
Moreover, the identity (\ref{10}) holds.
Then convergence (\ref{9}) follows from equality (\ref{1.11'}) and bound (\ref{2.4}).
\medskip\medskip

{\bf Acknowledgments}
This work was supported partly by the research grant of RFBR
(Grant No. 15-01-03587).
Author would thank Prof. A.I. Komech for helpful discussions.


\end{document}